\documentclass[bookmarksopen,twoside]{article}
\usepackage{graphicx,amssymb,mathrsfs,amsmath,color,fancyhdr}
\usepackage[bookmarks=true,dvipdfm]{hyperref}
\input vatola.sty \input cyracc.def
\input psfig.sty

\TagsOnRight

\renewcommand{\baselinestretch}{1.06} 
\catcode`@=11 \long\def\@makefntext#1{\noindent #1}
\newskip\tabcentering \tabcentering=1000pt plus 1000pt minus 1000pt
\def\REF#1{\par\hangindent\parindent\indent\llap{#1\enspace}}
\def\MCH#1#2{\setbox0=\hbox{\raise#1\hbox{#2}}\smash{\box0}}

\def\dl{\displaystyle}
\let\@oddfoot\@empty  \let\@evenfoot\@empty

\def\@evenhead{}\def\@oddhead{}
\def\@evenhead{\vbox{\hbox to \textwidth{\footnotesize\rm\hbox to
1.0cm{\thepage\hfill} \hfill\hspace{2mm}\footnotesize{
\emph{ZHANG Q. F. and CUI J. Z.}}}}}
\def\@oddhead{\vbox{\hbox to \textwidth{\footnotesize
{\it Regularity of the correctors and local
gradient estimate of the  homogenization} \hfill{\ } \hfill\hbox to
1cm{\hfill\thepage}}}}

\def\sec#1{\vspace{2mm}\noindent{{\bf #1}}\vspace{0.5mm}}
\def\th#1{\vspace{1mm}\noindent{\bf #1}\quad } 

\def\leq{\leqslant}

  \def\hml{\end{document}}  \newsymbol\wjzhml 203F \def\no{\noindent}

\begin{document}
\abovedisplayskip=3pt plus 1pt minus 1pt 
\belowdisplayskip=3pt plus 1pt minus 1pt 

\def\le{\leqslant}
\def\ge{\geqslant}
\def\dl{\displaystyle}




\vspace{8true mm}

\renewcommand{\baselinestretch}{1.9}\baselineskip 19pt
\begin{center}
\noindent{\LARGE\bf Regularity of the correctors and local
gradient estimate of the  homogenization for the elliptic equation: linear periodic case}

\vspace{0.5 true cm}

\noindent{\normalsize\sf ZHANG QiaoFu$^{\dag}$ \& CUI JunZhi
\footnotetext{\baselineskip 10pt
$^\dag$ Corresponding author\\
}}

\vspace{0.2 true cm}
\renewcommand{\baselinestretch}{1.5}\baselineskip 12pt
\noindent{\footnotesize\rm $
 Academy\, of\, Mathematics\, and\, Systems\, Science, \,
Chinese\, Academy\, of \,Sciences, \,Beijing\, 100190,\,China
$ \\
(email: zhangqf@lsec.cc.ac.cn,\,cjz@lsec.cc.ac.cn)\vspace{4mm}}
\end{center}

\baselineskip 12pt \renewcommand{\baselinestretch}{1.18}
\noindent{{\bf Abstract}\small\hspace{2.8mm} 
$C^\alpha$ and $W^{1,\infty}$ estimates for the f\mbox{}irst-order and second-order correctors in the homogenization
are presented based on the translation invariant and Li-Vogelius's gradient estimate  for the second order linear elliptic equation with piecewise smooth coef\mbox{}f\mbox{}icients.
If the data are smooth enough, the error of the f\mbox{}irst-order expansion for piecewise smooth coef\mbox{}f\mbox{}icients
is locally $O(\varepsilon)$ in the H\"older norm;
it is locally $O(\varepsilon)$ in $W^{1,\infty}$ when coef\mbox{}f\mbox{}icients are Lipschitz continuous.  It can be partly extended to the nonlinear parabolic equation.
}

\vspace{1mm} \no{\footnotesize{\bf Keywords:\hspace{2mm}}
gradient estimate, homogenization,
translation invariant,
De Giorgi-Nash  estimate
}

\no{\footnotesize{\bf MSC(2000):\hspace{2mm} } 35B27,\,35J65
 \vspace{2mm}
\baselineskip 15pt
\renewcommand{\baselinestretch}{1.22}
\parindent=10.8pt  
\rm\normalsize\rm

\sec{1\quad Introduction}
\newcommand{\dif}{\,\mathrm{d}}
\setcounter{section}{1}
\renewcommand\theequation{\arabic{section}.\arabic{equation}}

\newcounter{bh}
\usecounter{bh}
\newenvironment{remark}{\arabic{section}.\refstepcounter{bh}\arabic{bh}.\ }{}
\vspace{2mm}
\noindent Consider the homogenization of the following elliptic problem: f\mbox{}ind
$u_\varepsilon\in H^1_0(\Omega)$,
\begin{equation}\label{eq:fangcheng}
\mathcal{A}_\varepsilon u_\varepsilon\equiv
-\frac{\partial}{\partial x_i}\left(a_{ij}(\frac x \varepsilon)\frac{\partial u_\varepsilon}{\partial
x_j}\right)=f(x),\quad \mbox{in\,}\,\Omega\,.
\end{equation}
Here $\Omega\subset \mathbb{R}^n$ is  a bounded Lipschitz domain
and the summation convention is used. $A=(a_{ij})$ is symmetric and positive def\mbox{}inite; $a_{ij}(y)$ is
1-periodic in $y$, $1\leq i, j \leq n$;  $a_{ij}(y)$ is at least piecewise $C^\mu$ to obtain the error estimate in $C^\beta,W^{1,\infty}$.

Assume all of the data are smooth enough,
the $O(\varepsilon)$ error estimate in $L^\infty$  was presented by  A. Bensoussan, J. L. Lions and G. Papanicolaou  [1]; also see  M. Avellaneda and Lin FangHua [2].
O. A. Oleinik, A. S. Shamaev and G. A. Yosif\mbox{}ian  [3] proved the $O(\varepsilon^{1/2})$ estimate in $H^1$.
 Cao and Cui [4] studied the spectral properties and the numerical algorithms in perforated domains.
  Su \textit{et al} [5] investigated the quasi-periodic problems; Zhang and Cui [6]  gave a numerical example
  for the Rosseland equation.
  All of these were based on the multiple-scale expansion method [7].
There are also some other famous methods, such as periodic unfolding method [8],
 Multiscale Finite Element Method(MFEM [9]) and Heterogeneous Multiscale Method(HMM [10]).

The second-order expansion in Section 2 is classical which can be found in Chapter 1 [1] or Chapter 7 [7].
 Translation invariant  in Section 3 implies  the equivalence between the boundary and the interior estimate for an abstract periodic problem.
 The $C^\alpha,W^{1,\infty}$ estimates for correctors  in Section 4 follow from  the De Giorgi-Nash estimate and Li Yanyan-M. Vogelius's work
 for piecewise smooth   coef\mbox{}f\mbox{}icients, respectively. In Section 5, more than the traditional $L^\infty$ estimate ($a_{ij}(y)\in C^\gamma([0,1]^n)$, [2]), we obtain the $C^\beta$ error estimate
 ($a_{ij}(y)$ piecewise $ C^\mu$ in $C^{1,\alpha}$ subdomains, Corollary 5.4). At the end, we prove the main result: the error of the f\mbox{}irst-order expansion
is $O(\varepsilon)$ in $W^{1,\infty}_{loc}$ for Lipschitz continuous coef\mbox{}f\mbox{}icients (Corollary 6.3) based on M. Avellaneda-Lin FangHua's  gradient estimate. As far as we know, there are not such kinds ($C^\beta,W^{1,\infty}$) of error estimates  in the homogenization.

\vspace{5mm}

\sec{2\quad Second-order two-scale expansion}
\setcounter{equation}{0}
\setcounter{section}{2}
\setcounter{bh}{0}
\vspace{2mm}

\th{Def\mbox{}inition \remark}\textit{The periodic cell $Y=(0,1)^n$. Let $C^\infty_{per}(Y)$ be the subset of $C^\infty(\mathbb{R}^n)$ of Y-periodic functions(restricted on Y).
 Denote by $H^1_{per}(Y)$ the closure of $C^\infty_{per}(Y)$ in the $H^1$ norm.
$W^1_{per}(Y)=\{\varphi\in H^1_{per}(Y):\int_Y\varphi=0\}$. $\|u\|_{W^1_{per}(Y)}\equiv\|\nabla u\|_{L^2(Y)}$.
In the same way, we can def\mbox{}ine $W^1_{per}(Y_z)$ where $Y_z =Y+z,z\in \mathbb{R}^n$. }

 If $ u\in H^1_{per}(Y)$, $u$ has the same trace on the opposite faces of Y.
We look for a formal asymptotic expansion of the form
\begin{equation}\label{eq:ueps}
u_{\varepsilon}(x)=u_0(x) + \varepsilon u_1(x,\frac{x}{\varepsilon}) +
{\varepsilon}^2 u_2(x,\frac{x}{\varepsilon}) + ...
\end{equation}
where $u_1(\cdot,y),$ $u_2(\cdot, y)$ are Y-periodic
in y. Let $y=\frac{x}{\varepsilon}$, then
\begin{equation}\label{}
\frac{\partial}{\partial x_i} \rightarrow
\frac{\partial}{\partial x_i} + \frac{1}{\varepsilon} \frac{\partial}{\partial y_i}.
\end{equation}
 Substituting \eqref{eq:ueps} into \eqref{eq:fangcheng} and
 equating the power-like terms of $\varepsilon$, we introduce $N_m(y)\in W^1_{per}(Y)$, $1\leq m\leq n$, to make the terms of
 order $\varepsilon^{-1}$ equal zero,
\begin{equation}\label{eq:n}
\int_Y a_{ij}(y)\frac{\partial N_m}{\partial y_i}\frac{\partial \varphi}{\partial y_j} =
-\int_Y a_{mj}(y)\frac{\partial \varphi}{\partial y_j},\quad\forall \varphi(y) \in W^1_{per}(Y).
\end{equation}
Then let $u_1=N_m\partial_m u_0$. The problem for $u_2$ (the part of order
$\varepsilon^{0}$) admits a unique solution if and only if there exists a $u_0\in  H^1_0(\Omega)$ such that(a compatibility condition, see Theorem 4.26 [7])
\begin{equation}\label{eq:u0}
-\frac{\partial}{\partial x_i}
 [a_{ij}^0
\frac{\partial u_0}{\partial x_j}]
 =
f,\quad a_{ij}^0 =\int_Y[a_{ij}(y)+a_{il}(y)\frac{\partial N_j}{\partial y_l}]\mbox{d} y.
\end{equation}
This equation called the homogenization equation
is well-posed because $(a^0_{ij})$ is elliptic (Proposition 6.12 [7]).

Find  $M_{kl}\in W^1_{per}(Y),1\leq k,l\leq n,$
 such that
\begin{equation}\label{eq:m}
\int_Y
a_{ij}(y)\frac{\partial M_{kl}}{\partial y_i}\frac{\partial \varphi}{\partial y_j} = \int_Y
\left[a_{kl}+a_{km}\frac{\partial N_l}{\partial y_m} -a_{kl}^0\right]\varphi -
 \int_Ya_{ik}N_l\frac{\partial \varphi}{\partial y_i},\quad\forall  \varphi\in W^1_{per}(Y).
\end{equation}

If $\varphi(y)=1$, the righthand side of the above equation equals zero (a compatibility condition). So
 let
$u_2=M_{kl}\partial^2_{kl} u_0 $ to make the  $O(\varepsilon^{0})$ terms equal zero. Note that $N_m,M_{kl},u_0$ are independent of $\varepsilon$. We
 will use this fact again and again.

....

\th{Corollary \remark}
\textit{Under the hypotheses of Theorem 6.2,}
\begin{equation}\label{}
\sup_{\overline{\Omega'}}|\nabla (u_\varepsilon-u_0-\varepsilon u_1)|\leq
 C\varepsilon,\quad \Omega'\subset\subset \Omega;
\end{equation}
\begin{equation}\label{eq:flux}
\sup_{\overline{\Omega'}}|A(\frac x \varepsilon)\nabla (u_\varepsilon-u_0-\varepsilon u_1)|\leq
 C\varepsilon,\quad \Omega'\subset\subset \Omega.
\end{equation}
\th{Proof.} (1) We only need to prove that $|\varepsilon\partial_i u_2|=|\partial_i(\varepsilon M_{kl}(\frac x \varepsilon)\partial_{kl}^2u_0)|\leq C$.
\begin{equation}\label{}
 \varepsilon\frac{\partial u_2}{\partial x_i}=\frac{\partial M_{kl}}{\partial y_i}(\frac x \varepsilon)\partial_{kl}^2u_0+\varepsilon  M_{kl}(\frac x \varepsilon)\partial_{ikl}^3u_0.
\end{equation}
$M_{kl},\frac{\partial M_{kl}}{\partial y_i}$ are bounded from Theorem 4.4; $\partial_{ikl}^3u_0\in W^{1,q}(\Omega')\hookrightarrow C^{0,\alpha}(\overline{\Omega'})\subset L^\infty(\Omega')$.

(2)Note that $A(\frac x \varepsilon)=(a_{ij})$ is bounded.
\hfill   $\square$

\th{Remark \remark}We give the estimate \eqref{eq:flux} because the f\mbox{}lux ($A\nabla u$) is very important in physics.
One can consider the tensor case where the f\mbox{}lux may be the stress in linear elasticity.

\vspace{5mm}

 \sec{7\quad Some  problems}
\setcounter{equation}{0}
\setcounter{section}{7}
\setcounter{bh}{0}

\vspace{2mm}
 It is possible to partly extend  the above results to the following cases:

(1) tensor case: Avellaneda-Lin's Lemma 6.1 is true for the tensor case [2] and elliptic systems with Neumann boundary
 conditions [13]; Li-Vogelius's work was extended  in [14].

(2) nonlinear case:   Fusco and Moscariello [15] studied the homogenization of quasilinear divergence structure operators. For the second-order expansion, see [16].

 (3) parabolic case: the parabolic $C^{\alpha,\alpha/2}$ estimate under mixed boundary conditions was presented in [17];
 Li-Vogelius's gradient estimate was extended to parabolic systems  in [14].

(4) nonsmooth case: if the domain is only convex or the righthand side is piecewise smooth,
 there are many interesting problems. One problem is that the hypotheses in Theorem 5.3 are very strong: $\partial \Omega\in C^{2,1}$, $f\in W^{1,q}(\Omega)$, $q>n$.
This is a common dif\mbox{}f\mbox{}iculty for the multiple-scale method (see  [1], [2]).

(5) How can we get a global $W^{1,\infty}$ error estimate with a proper boundary
 corrector?


Some results will appear elsewhere.

\vspace{3mm}\th{Acknowledgements}
This work is supported by National Natural Science Foundation of China (Grant No. 90916027).
The authors thank Professor Yan NingNing and the referees for their careful reading and helpful comments.

\vskip0.1in \no {\normalsize \bf References}
\vskip0.1in\parskip=0mm \baselineskip 15pt
\renewcommand{\baselinestretch}{1.15}

\footnotesize\parindent=6mm
 \REF{1\ }Bensoussan A., Lions J. L., Papanicolaou G.,
Asymptotic Analysis for Periodic Structures, Amsterdam: North-Holland, 1978

\footnotesize\parindent=6mm
 \REF{2\ }Avellaneda M., Lin F. H., Compactness method in the theory of homogenization,
Comm Pure Appl Math, 1987, 40(6): 803-847

\footnotesize\parindent=6mm
 \REF{3\ }Oleinik O. A., Shamaev A. S., Yosifian G. A.,
Mathematical Problems in Elasticity and Homogenization, Amsterdam: North-Holland, 1992

\footnotesize\parindent=6mm
 \REF{4\ }Cao L. Q., Cui J. Z., Asymptotic expansions and numerical algorithms of eigenvalues and eigenfunctions of the Dirichlet problem for second order elliptic equations in
perforated domains, Numer Math, 2004, 96(3): 525-581

\footnotesize
\parindent=6mm
 \REF{5\ }Su F., Cui J. Z., Xu Z., \textit{et al},
A second-order and two-scale computation method for the quasi-periodic
structures of composite materials,
Finite Elements in Analysis and Design, 2010,
46(4): 320-327

\footnotesize\parindent=6mm
 \REF{6\ }Zhang Q. F., Cui J. Z.,
Multi-scale analysis method for combined conduction-radiation heat transfer of periodic composites,
 Advances in Heterogeneous Material Mechanics(eds. Fan J. H., Zhang J. Q., Chen H. B., \textit{et al}), Lancaster: DEStech Publications, 2011, 461-464

\footnotesize\parindent=6mm
 \REF{7\ }Cioranescu D., Donato P.,     An Introduction to Homogenization, Oxford: Oxford    University Press, 1999

\footnotesize\parindent=6mm
 \REF{8\ }Griso G., Error estimate and unfolding for periodic
homogenization, Asymptotic Analysis, 2004, 40(3): 269-286

\footnotesize
\parindent=6mm
 \REF{9\ }Hou T. Y., Wu X. H.,
A multiscale finite element method for elliptic problems in composite materials and porous media,
 J Comput Phys, 1997, 134(1): 169-189

\footnotesize
\parindent=6mm
 \REF{10\ }E W. N., Engquist B., Huang Z. Y.,
Heterogeneous multiscale method: a general methodology for multiscale modeling,
Phys Rev B, 2003, 67(092101): 1-4

\footnotesize\parindent=6mm
 \REF{11\ }Gilbarg D., Trudinger N. S., Elliptic Partial Differential Equations  of Second Order, Berlin: Springer, 2001

 \footnotesize\parindent=6mm
 \REF{12\ }Li Y. Y., Vogelius M., Gradient estimates for solutions to divergence form elliptic equations with
    discontinuous coef\mbox{}f\mbox{}icients, Arch Rational Mech Anal, 2000, 153(2): 91-151.

 \footnotesize\parindent=6mm
 \REF{13\ }Kenig C. E., Lin F. H., Shen Z. W., Homogenization of elliptic systems with Neumann boundary
    conditions, 2010,  arXiv: 1010.6114v1, [math.AP]

\footnotesize\parindent=6mm
 \REF{14\ }Li H. G., Li Y. Y., Gradient estimates for parabolic systems from
composite material, 2011,
arXiv: 1105.1437v1,  [math.AP]

\footnotesize\parindent=6mm
 \REF{15\ }Fusco N.,  Moscariello G., On the homogenization of
quasilinear divergence structure operators,
Annali di Matematica Pura ed Applicata, 1986, 146(1): 1-13

\footnotesize\parindent=6mm
 \REF{16\ }Zhang Q. F., Cui J. Z.,  Error estimate of the second-order homogenization for divergence-type nonlinear elliptic equation,
2011,
arXiv: 1108.5070v1, [math-ph]

\footnotesize\parindent=6mm
 \REF{17\ }Griepentrog  J. A., Recke L., Local existence, uniqueness and smooth dependence for
     nonsmooth quasilinear parabolic problems, J Evol Equ, 2010, 10(2): 341-375

 \hml